\documentclass[12pt]{article}
\usepackage[final]{epsfig}	
\usepackage{amsfonts,amssymb,latexsym,amscd, theorem,epsfig,graphics}

\newtheorem{theorem}{Theorem}
\newtheorem{lemma}{Lemma}[theorem]

\newtheorem{example}{Example}
\newtheorem{proposition}{Proposition}[example]

\begin{document}
\newcommand{\eps}{{\varepsilon}}
\newcommand{\proofend}{$\Box$\bigskip}
\newcommand{\C}{{\mathbf C}}
\newcommand{\Q}{{\mathbf Q}}
\newcommand{\R}{{\mathbf R}}
\newcommand{\Z}{{\mathbf Z}}
\newcommand{\RP}{{\mathbf {RP}}}

\title {Non-existence of $n$-dimensional $T$-embedded discs in ${\bf R}^{2n}$}
\author{Gordana Stojanovic and Serge Tabachnikov \thanks{Partially supported by NSF}\\
{\it Department of Mathematics, Brown University}\\
{\it Providence, RI 02912, USA}\\
and\\
{\it Department of Mathematics, Penn State University}\\
{\it University Park, PA 16802, USA}
}
\date{}
\maketitle

A number of recent papers concerned various non-degeneracy conditions on embedding and immersions of smooth manifolds in affine and projective spaces defined in terms of mutual positions of the tangent spaces at distinct points, see \cite{Gh1, Gh2, G-S, G-T, S-S, Ta, Wu}. Following Ghomi \cite{Gh1}, a $C^1$-embedded manifold $M^n \subset \R^N$ is called  $T$-embedded if the tangent spaces to $M$ at distinct points do not intersect. For example, the cubic curve $(x,x^2,x^3)$ is a $T$-embedding of $\R$ to $\R^3$, and the direct product of such curves gives a $T$-embedding of $\R^n$ to $\R^{3n}$. 

A $T$-embedding $M^n \to \R^N$ induces a topological embedding of the tangent bundle $TM \to \R^N$, hence $N\geq 2n$.
One of the results in \cite{Gh1} is that  no closed manifold $M^n$  admits $T$-embeddings to $\R^{2n}$. In this note we strengthen this result as follows. 

\begin{theorem} \label{main}
There exist no $C^2$-smooth $T$-embedded discs $D^n$ in $\R^{2n}$.
\end{theorem}

\noindent {\bf Proof}. Arguing by contradiction, assume that such a disc $D^n$ exists. Choose the tangent space at the origin and its orthogonal complement as coordinate $n$-dimensional spaces. Making $D$ smaller, if necessary, assume that the disc is the graph of a (germ of a) $C^2$ smooth map $f: {\bf R}^n \to {\bf R}^n$. Let $U\subset {\bf R}^n$ be the domain of $f$.

Let $z=(u,f(u)) \in D$ where $u\in U$. The tangent space $T_z D$ is given by a linear equation $y=A(u) x - b(u)$ where $A(u)$ is an $n\times n$ matrix and $b(u)$ is a vector in ${\bf R}^n$, both depending on $u$. In terms of $f$, they have the following expressions. Let $f_1,...,f_n$ be the components of $f$.

\begin{lemma} \label{L1} One has: 
$$
A_{ij}=\frac{\partial f_i}{\partial u_j}, \quad  b_i=\sum_{k=1}^n \frac{\partial f_i}{\partial u_k} u_k - f_i.
$$ 
\end{lemma}

\noindent {\bf Proof}. The first statement is obvious, and the second follows from the fact that the space  $y=A(u) x - b(u)$ passes through the point $z=(u,f(u))$. \proofend

One has the next characterization of  $t$-discs.

\begin{lemma} \label{L2} For all $u\neq v \in U$, the vector $b(u)-b(v)$ does not belong to ${\rm Im}(A(u)-A(v))$. 
\end{lemma}

\noindent {\bf Proof}. The spaces  $y=A(u) x - b(u)$ and $y=A(v) x - b(v)$ intersect if and only if
$b(u)-b(v) \in {\rm Im}(A(u)-A(v))$. \proofend

\begin{lemma} \label{L3} 
If $u\neq v$ then $b(u)\neq b(v)$ and $A(u)-A(v)$ is degenerate.
\end{lemma}

\noindent {\bf Proof}. The first claim follows from the fact that zero vector lies in any subspace, contradicting Lemma \ref{L2}. If $A(u)-A(v)$ is nondegenerate then it is surjective, again contradicting Lemma \ref{L2}.
\proofend

Now we compute the Jacobian of the map $b:U \to \R^n$. Denote by $E$ the Euler vector field in $\R^n$:
$$
E=\sum_{k=1}^n u_k \frac{\partial}{\partial u_k}.
$$

\begin{lemma} \label{L4}
 One has: 
 $$
 \frac{\partial b_i}{\partial u_j} = \sum_k \frac{\partial^2 f_i}{\partial u_j \partial u_k} u_k=E(A_{ij}).
 $$
 \end{lemma}

\noindent {\bf Proof}. This follows from Lemma 1.  \proofend

\begin{lemma} \label{L5} 
For all $u\in U$, the Jacobian $Jb$ of the map $b$ is degenerate.
\end{lemma}

\noindent {\bf Proof}.  Lemma \ref{L4} implies that 
$$
Jb=\lim_{\varepsilon \to 0} \frac{A(u+\varepsilon u)-A(u)}{\varepsilon}.
$$
By Lemma \ref{L3} with $v=u+\varepsilon u$, the numerator is a degenerate matrix for all $\varepsilon$, and so is its quotient by $\varepsilon$. Thus $Jb$ is a limit of degenerate matrices. Since determinant is a continuous function, the limit  also has zero determinant and therefore is degenerate.  \proofend

Finally, we arrive at a contradiction. By Lemma \ref{L3}, the map $b$ is one-to-one, and by the invariance of domain theorem, its image has positive measure. By Lemma \ref{L5}, every value of $b$ is singular, and by Sard's Lemma its image has zero measure. This completes the proof of Theorem \ref{main}.
\proofend

 According to Lemma \ref{L3}, the $n$-parameter family of $n\times n$ matrices $A(u), u\in D^n$ enjoys the property that $A(u)-A(v)$ is degenerate for all $u\neq v$. If $n=2$, such families can be explicitly described. Assume that not all matrices  $A(u)$ are zero.

\begin{theorem} \label{family}
The family $A(u)$ consists either of the matrices with a fixed $1$-dimensional image or with a fixed $1$-dimensional kernel.
\end{theorem}

\noindent {\bf Proof}. Let $M_2$ be the space of linear maps $\R^2 \to \R^2$. One has a non-degenerate quadratic form in $M_2$ given by the determinant of a matrix; this form has signature $(2,2)$. Consider the respective dot product.

Let $V \subset M_2$ be the linear span of the family $A(u)$. 

\begin{lemma} \label{isotr}
The subspace $V$ is isotropic.
\end{lemma}

\noindent {\bf Proof}. It suffices to prove that  $A(u) \cdot A(v)=0$ for all $u,v$. If $u=v$, this means precisely that $A(u)$ is degenerate.  For $u \neq v$, the matrix $A(u)-A(v)$ is degenerate, hence $(A(u)-A(v)) \cdot (A(u)-A(v)) =0$. Using bilinearity of the dot product, it follows that $A(u) \cdot A(v) =0$.
\proofend

Since the dot product is non-degenerate, an isotropic subspace is at most $2$-dimensional. 

\begin{lemma} \label{class}
A 2-dimensional isotropic subspace in $M_2$ consists either of the matrices with a fixed $1$-dimensional image or with a fixed $1$-dimensional kernel.
\end{lemma}

\noindent {\bf Proof}. 
Let $A \in V$ be a non-zero matrix. Choose a basis in the target space $\R^2$ in such a way that Im $A$ is orthogonal to the column vector $(0,1)$. Then 
$$
 A=\left(\begin{array}{cc}
 a&b\\
 0&0
 \end{array}\right)
$$
with $a^2+b^2\neq 0$. Let $B\in V$ be another matrix, not proportional to $A$. Then $A \cdot B=0$, and hence $$
 B=\left(\begin{array}{cc}
 c&d\\
 at&bt
 \end{array}\right)
$$
for some real $c,d,t$. If $t=0$ then $(c,d)$ is not proportional to $(a,b)$, and the space $V$ consists of matrices with zero second row. This is the first case of the lemma: the matrices have a fixed image spanned by the column vector $(1,0)$.

Otherwise, $t\neq 0$. Since $B$ is degenerate, one has: $(c,d)=s(a,b)$ for some real $s$. Then 
$$
\frac{B-sA}{t}= \left(\begin{array}{cc}
 0&0\\
 a&b
 \end{array}\right),
$$
and the space $V$ consists of matrices with a fixed kernel spanned by the column vector $(-b,a)$.
\proofend

Lemma \ref{class} obviously implies Theorem \ref{family}.
\proofend

For $n=2$, Theorem \ref{family} implies the claim of Theorem \ref{main}. Indeed, assume that the Jacobi matrix $Jf$ has a fixed $1$-dimensional kernel, say, spanned by vector $\xi$. Then the map $f$ has zero directional derivative along $\xi$, and the tangent planes to the graph of $f$ are the same along this direction. Hence this graph is not $T$-embedded. Likewise, if  $Jf$ has a fixed $1$-dimensional image then the transpose matrix  has a fixed kernel, say, $\eta$. This implies that the function $f(u) \cdot \eta$ has zero differential, and hence the  image of $f$ is $1$-dimensional. It follows that the graph of $f$ belongs to a $3$-dimensional space and therefore is not $T$-embedded.

Let us conclude with two examples motivated by the following erroneous attempt to prove Theorem \ref{main}: if there exists a $T$-embedded disc $D^n \subset \R^{2n}$ then 
its tangent spaces provide a foliation $\cal F$ of a domain in $\R^{2n}$ by $n$-dimensional affine subspaces. Then $D^n$ is everywhere tangent to the leaves of this $n$-dimensional foliation  and therefore  must lie within a leaf. The mistake in this argument is that, no matter how smooth the embedding is, the foliation $\cal F$ is not differentiable. This phenomenon is illustrated in the following example.

\begin{example} \label{ex1}
{\rm Let $\gamma$ be a smooth plane curve with positive curvature and free from vertices (extrema of curvature). Then, by the classical Kneser theorem (1912), the osculating circles to $\gamma$ are pairwise disjoint and nested as illustrated in figure \ref{osc}; see, e.g.,  \cite{Gu,Ta1}. These osculating circles foliate the annulus $A$ between the largest and smallest of them.  Denote this foliation by $\cal F$. Then $\cal F$ is not $C^1$, namely, one has the following result.

\begin{figure}[ht]
\centerline{\epsfbox{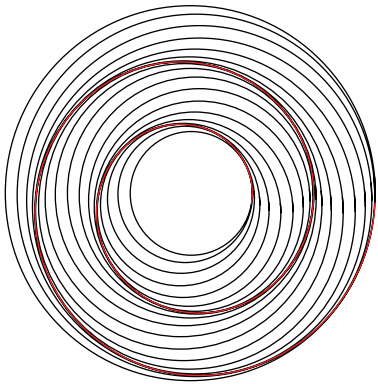}}
\caption{Osculating circles of a spiral}
\label{osc}
\end{figure}

\begin{proposition} \label{kneser}
Let $f: A \to \R$ be a differentiable function, constant on the leaves of $\cal F$. Then $f$ is constant in $A$. 
\end{proposition}

\noindent {\bf Proof}. Since $f$ is constant on the leaves of $\cal F$, the differential $df$ vanishes on any vector tangent to any leaf. Since $\gamma$ is everywhere tangent to the leaves, $df$ is zero on the tangent vectors to $\gamma$. Hence $f$ is constant on $\gamma$. But $A$ is the union of the leaves of $\cal F$ through the points of $\gamma$, hence $f$ is constant in $A$.  
\proofend
}
\end{example}

One  also wonders whether $\R^{2n}$ can be foliated by non-parallel affine $n$-dimensional subspaces (clearly impossible for $n=1$). 

\begin{example} \label{ex2}
{\rm The following construction gives a foliation of $\R^4$ by pairwise non-parallel $2$-dimensional affine subspaces.  Start with partitioning $3$-dimensional space  into the vertical $z$-axis and the hyperboloids of 1 sheet 
$$
x^2+y^2=t(z^2+1), \ \ t>0
$$ 
(when $t=0$, one has the $z$-axis). Each hyperbolid is foliated by lines, and thus $\R^3$ gets foliated by lines; these lines are pairwise skew. Multiply this foliation by $\R^1$ to obtain the desired example.
}
\end{example}

\end{document}